\documentclass[11pt,reqno]{amsart}
\usepackage{amsfonts, amssymb, amsmath, enumerate, bm, xspace, graphicx, caption, subcaption}
\usepackage{pbox}
\usepackage{epstopdf}

\numberwithin{equation}{section}

\DeclareMathOperator{\E}{\mathbb{E}}
\DeclareMathOperator{\diag}{diag}

\def \P {\mathbb{P}}
\def \R {\mathbb{R}}

\def \LL {\mathcal{L}}

\def \e {\varepsilon}

\def \s {\sigma}

\def \tran {\mathsf{T}}

\def \one {{\textbf 1}}

\newtheorem{theorem}{Theorem}[section]
\newtheorem{proposition}[theorem]{Proposition}
\newtheorem{corollary}[theorem]{Corollary}

\theoremstyle{remark}

\title[]{Concentration of random graphs and application to community detection}

\author{Can M. Le, Elizaveta Levina \and Roman Vershynin}

\address{Department of Statistics, University of California, Davis, One Shields Ave, Davis, CA 95616, U.S.A.}
\email{canle@ucdavis.edu}

\address{Department of Statistics, University of Michigan, 1085 S. University Ave, Ann Arbor, MI 48109, U.S.A.}
\email{elevina@umich.edu}

\address{Department of Mathematics, University of California, Irvine,   	 
340 Rowland Hall, Irvine, CA 92697, U.S.A.}
\email{rvershyn@uci.edu}

\thanks{}

\begin{document}


\begin{abstract}
 Random matrix theory has played an important role in recent work on statistical network analysis.   In this paper, we review recent results on regimes of concentration of random graphs around their expectation, showing that dense graphs concentrate and sparse graphs concentrate after regularization.     We also review relevant network models that may be of interest to probabilists considering directions for new random matrix theory developments, and random matrix theory tools that may be of interest to statisticians looking to prove properties of network algorithms.   Applications of concentration results to the problem of community detection in networks are discussed in detail.   
\end{abstract}


\maketitle


\section{Introduction}
A lot of recent interest in concentration of random graphs has been generated by problems in network analysis, a very active interdisciplinary research area with contributions from probability, statistics, physics, computer science, and the social sciences all playing a role.     Networks represent relationships (edges) between objects (nodes), and a network between $n$ nodes is typically represented by its $n \times n$ adjacency matrix $A$.    We will focus on the case of simple undirected networks, where $A_{ij} = 1$ when nodes $i$ and $j$ are connected by an edge, and 0 otherwise, which makes $A$ a symmetric matrix with binary entries.   It is customary to assume the graph contains no self-loops, that is, $A_{ii} = 0$ for all $i$, but this is not crucial.    In general, networks may be directed ($A$ is not symmetric), weighted (the entries of $A$ have  a numerical value representing the strength of connection), and/or signed (the entries of $A$ have a sign representing whether the relationship is positive or negative in some sense).

Viewing networks as random realizations from an underlying network model enables analysis and inference, with the added difficulty that we often only observe a single realization of a given network.    Quantities of interest to be inferred from this realization may include various network properties such as the node degree distribution, the network radius, and community structure.   Fundamental to these inferences is the question of how close a single realization of the matrix $A$ is to the population mean, or the true model, $\E A$.   If $A$ is close to $\E A$, that is, $A$ concentrates around its mean, then inferences drawn from $A$ can be transferred to the population with high probability.    

In this paper, we aim to answer the question ``When does $A$ concentrate around $\E A$?'' under a number of network models and asymptotic regimes.   We also show that in some cases when the network does not concentrate, a simple regularization step can restore concentration.     While the question of concentration is interesting in its own right, we especially focus on the implications for the problem of community detection, a problem that has attracted a lot of attention in the networks literature.   When concentration holds, in many cases a simple spectral algorithm can recover communities, and thus concentration is of practical and not only theoretical interest.



\section{Random network models}\label{sec:models}

Our concenrtation results hold for quite general models, but, for the sake of clarity, we provide a brief review of network models, starting from the simplest model and building up in complexity.  

\subsection*{The  Erd{\H{o}}s-R{\'e}nyi  (ER) graph}
The simplest random network model is the Erd{\H{o}}s-R{\'e}nyi graph $G(n,p)$ \cite{Erdos&Renyi1959}. Under this model, edges are independently drawn between all pairs of nodes according to a Bernoulli distribution with success probability $p$. Although the ER model provides an important building block in network modeling and is attractive to analyze, it almost never fits network data observed in practice. 

\subsection*{The stochastic block model (SBM)}    The SBM is perhaps the simplest network model with community structure, first proposed by \cite{Holland83}. Under this model, each node belongs to exactly one of $K$ communities,  and the node community membership $c_i$ is drawn independently from a multinomial distribution on $\{1, \dots, K\}$ with probabilities $\pi_1, \dots, \pi_K$. Conditional on the label vector $c$, edges are drawn independently between each pair of nodes $i,j$, with  
$$\P (A_{ij} = 1) = B_{c_ic_j},$$ 
where $B$ is a symmetric $K\times K$ matrix controlling edge probabilities. Note that each community within SBM is an ER graph.   The main question of interest in network analysis is estimating the label vector $c$ from $A$, although model parameters $\pi$ and $P$ may also be of interest.

\subsection*{The degree-corrected stochastic block model (DCSBM)}  While the SBM does incorporate community structure, the assumption that each block is an ER graph is too restrictive for many real-world networks.   In particular, ER graphs have a Poisson degree distribution, and real networks typically fit the power law or another heavy-tailed distribution better, since they often have ``hubs'', influential nodes with many connections.   An extension removing this limitation, the degree-corrected stochastic block model (DCSBM) was proposed by \cite{Karrer10}.   The DCSBM is like an SBM but with each node assigned an additional parameter $\theta_i > 0$ that controls its expected degree, and edges drawn independently with $$\P (A_{ij} = 1) = \theta_i \theta_j B_{c_ic_j}.$$ 
Additional constraints need to be imposed on $\theta_i$ for model identifiability;  see \cite{Karrer10,Zhaoetal2012} for options. 

\subsection*{The latent space model (LSM)}  Node labels under the SBM or the DCSBM can be thought of as latent (unobserved) node positions in a discrete space of $K$ elements.  More generally, latent positions can be modeled as coordinates in $\R^d$, or another set equipped with a distance measure.   The LSM \cite{Hoff2002} assumes that each node $i$ is associated with an unknown position $x_i$ and edges are drawn independently between each pair of nodes $i,j$ with probability inversely proportional to the distance between $x_i$ and $x_j$. If latent positions $x_i$ form clusters (for example, if they are drawn from a mixture of Gaussians), then a random network generated from this model exhibits community structure.   Inferring the latent positions can in principle lead to insights into how the network was formed, beyond simple community assignments.  

\subsection*{Exchangeable random networks} 
An analogue of de Finetti's theorem for networks, due to Hoover and Aldous \cite{Hoover1979, Aldous81}, shows that any network whose distribution is invariant under node permutations can be represented by 
$$
A_{ij} = g(\alpha,\xi_i,\xi_j,\lambda_{ij}),
$$ 
where $\alpha$, $\xi_i$ and $\xi_j$ are independent and uniformly distributed on $[0,1]$, and $g(u,v,w,z)=g(u,w,v,z)$ for all $u,v,w,z$.    This model covers all the previously discussed models as special cases, and the function $g$, called the graphon, can be estimated up to a permutation under additional assumptions;  see   
\cite{Olhede2013,Gao&Lu&Zhou2015,Zhang&Levina&Zhu2017}.

\subsection*{Network models with overlapping communities} In practice, it is often more reasonable to allow nodes to belong to more than one community.   Multiple such models have been proposed, including the Mixed Membership Stochastic Block Model (MMSBM) \cite{Airoldi2008}, the Ball-Karrer-Newman Model (BKN) \cite{Ball&Karrer&Newman2011}, and the OCCAM model \cite{Zhang&Levina&Zhu2014}.   MMSBM allows different memberships depending on which node the given node interacts with;  the BKN models edges as a sum of multiple edges corresponding to different communities;  and OCCAM relaxes the membership vector $c$ under the SBM to have entries between 0 and 1 instead of exactly one ``1''.   All of these models are also covered by the theory we present, because, conditional on node memberships, all these networks are distributed according to an inhomogeneous  Erd{\H{o}}s-R{\'e}nyi model, the most general model we consider,  described next.




\subsection*{The inhomogeneous Erd{\H{o}}s-R{\'e}nyi model}
All models described above share an important property: conditioned on node latent positions, edges are formed independently.    The most general form of such a model is the inhomogeneous Erd{\H{o}}s-R{\'e}nyi model (IERM) \cite{Bollobas2007}, where each edge is independently drawn, with $\P(A_{ij} = 1) = P_{ij}$, where $P=(P_{ij}) = \E A$.   Evidently, additional assumptions have to be made if latent positions of nodes (however they are defined) are to be recovered from a single realization of $A$.   We will state concentration results under the IERM as generally as possible, and then discuss additional assumptions under which latent positions can also be estimated reliably.

\subsection*{Scaling}
We have so far defined all the models for a fixed number of nodes $n$, but in order to talk about concentration, we need to determine how the expectation $P_n = \E A_n$ changes with $n$.   Most of the literature defines 
$$P_n = \rho_n P$$
where $P$ is a matrix with constant non-negative entries,  and $\rho_n$ controls the average expected degree of the network, $d = d_n = n\rho_n$.    Different asymptotic regimes have been studied, especially under the SBM; see \cite{AbbeReview2017} for a review.   Unless $\rho_n \rightarrow 0$, the average network degree $d = \Omega(n)$, and the network becomes dense as $n$ grows.    In the SBM literature, the regime $d_n \gg \log n$ is sometimes referred to as semi-dense;   $d_n \rightarrow \infty$ but not faster than $\log n$ is semi-sparse;  and the constant degree regime $d_n = O(1)$ is called sparse.   We will elaborate on these regimes and their implications later on in the paper.

%

\section{Useful random matrix results}

We start from presenting a few powerful and general tools in random matrix theory which can help prove concentration bounds for random graphs.

\begin{theorem}[Bai-Yin law \cite{Bai&Yin1988}; see \cite{FurKom80} an for earlier result]
  Let $M = (M_{ij})_{i,j=1}^\infty$ be an infinite, symmetric, and diagonal-free random matrix 
  whose entries above the diagonal are i.i.d. random variables with zero mean and variance $\s^2$.
  Suppose further that $\E M_{ij}^4 < \infty$. 
  Let $M_n = (M_{ij})_{i,j=1}^n$ denote the principal minors of $M$. Then, as $n \to \infty$, 
  \begin{equation}		\label{eq: 2}
  \frac{1}{\sqrt{n}} \|M_n\| \to 2 \quad \text{almost surely.}
  \end{equation}
\end{theorem}


\begin{theorem}[Matrix Bernstein's inequality]  \label{thm: matrix Bernstein}
  Let $X_1,\ldots,X_N$ be independent, mean zero, $n \times n$ symmetric random matrices, 
  such that $\|X_i\| \le K$ almost surely for all $i$.
  Then, for every $t \ge 0$ we have
  $$
  \P \Big\{ \Big\| \sum_{i=1}^N X_i \Big\| \ge t \Big\} 
  \le 2n \exp \Big( -\frac{t^2/2}{\s^2 + Kt/3} \Big).
  $$
  Here $\s^2 = \left\| \sum_{i=1}^N \E X_i^2 \right\|$ is the norm of the ``matrix variance'' of the sum.
\end{theorem}

\begin{corollary}[Expected norm of sum of random matrices]  
\label{cor: matrix Bernstein expectation}
  We have  
  $$
  \E \Big\| \sum_{i=1}^N X_i \Big\| 
  \lesssim \s \sqrt{\log n} + K \log n.
  $$
\end{corollary}

The following result gives sharper bounds on random matrices than matrix Bernstein's inequality,   but requires independence of entries.   
\begin{theorem}[Bandeira-van Handel \cite{Bandeira&Handel2016} Corollary~3.6]
\label{Bandeira}
  Let $M$ be an $n \times n$ symmetric random matrix with independent entries 
  on and above the diagonal. Then 
  $$
  \E \|M\| \lesssim \max_i \Big( \sum_j \s_{ij}^2 \Big)^{1/2} + \sqrt{\log n} \, \max_{i,j} K_{ij},
  $$
  where $\s_{ij}^2 = \E M_{ij}^2$ are the variances of entries and $K_{ij} = \|M_{ij}\|_\infty$.
\end{theorem}

\begin{theorem}[Seginer's theorem \cite{Seginer2000}]		\label{thm: Seginer}
  Let $M$ be a $n \times n$ symmetric random matrix 
  with i.i.d.\ mean zero entries above the diagonal
  and arbitrary entries on the diagonal. 
  Then 
  $$
  \E \|M\| \asymp \E \max_i \|M_i\|_2
  $$
  where $M_i$ denote the columns of $M$.
\end{theorem}

The lower bound in Seginer's theorem is trivial; it follows from the fact that the operator norm of a matrix
is always bounded below by the Euclidean norm of any of its columns. 
The original paper of Seginer \cite{Seginer2000} proved the upper bound for non-symmetric matrices with independent entries. The present statement of Theorem~\ref{thm: Seginer}
can be derived by a simple symmetrization argument, 
see \cite[Section~4.1]{Hajek&Wu&Xu2016}. 

\section{Dense networks concentrate}
If $A = A_n$ is the adjacency matrix of a $G(n,p)$ random graph with a constant $p$, 
then the Bai-Yin law gives
$$
\frac{1}{\sqrt{n}} \|A - \E A\| \to 2 \sqrt{p(1-p)}.
$$ 
In particular, we have
\begin{equation}		\label{eq: Bai-Yin}
\|A - \E A\| \le 2 \sqrt{d}
\end{equation}
with probability tending to one, where $d = np$
is the expected node degree.

Can we expect a similar concentration for sparser Erd\"os-R\'enyi graphs, 
where $p$ is allowed to decrease with $n$? 
The method of 
\cite{Friedman&Kahn&Szemeredi1989} adapted by Feige and Ofek \cite{FeiOfe05} gives
\begin{equation}		\label{eq: sqrt d}
\|A - \E A\| = O(\sqrt{d})
\end{equation}
under the weaker condition $d \gtrsim \log n$, which is optimal, as we will see shortly.
This argument actually yields \eqref{eq: sqrt d} for inhomogeneous random graphs $G(n, (p_{ij}))$ as well, 
and for $d = \max_{ij} np_{ij}$, see e.g. \cite{Lei&Rinaldo2015, Chin&Rao&Vu2015}. 

Under a weaker assumption $d = np \gg \log^4 n$, 
Vu \cite{Vu2007} proved a sharper bound for $G(n,p)$, namely 
\begin{equation}		\label{eq: optimal concentration}
\|A - \E A\| = (2+o(1)) \sqrt{d},
\end{equation}
which essentially extends \eqref{eq: Bai-Yin} to sparse random graphs.
Very recently, Benaych-Georges, Bordenave and Knowles \cite{Benaych-georges&Bordenave&Knowles2017} 
were able to derive \eqref{eq: optimal concentration} under the optimal condition $d \gg \log n$. 
More precisely, they showed that if $4 \le d \le n^{2/13}$, then 
$$
\E \|A - \E A\| \le 2 \sqrt{d} + C \sqrt{ \frac{\log n}{1 + \log(\log(n)/d)}}.
$$
The argument of \cite{Benaych-georges&Bordenave&Knowles2017} applies more generally to inhomogeneous random graphs
$G(n,(p_{ij}))$ under a regularity condition on the connection probabilities $(p_{ij})$. 
It even holds for more general random matrices that may not necessarily have binary entries.


%

To apply Corollary \ref{cor: matrix Bernstein expectation} to the adjacency matrix $A$ of an ER random graph $G(n,p)$, decompose $A$ into a sum of independent random matrices
$A = \sum_{i \le j} X_{ij}$,
where each matrix $X_{ij}$ contains a pair of symmetric entries of $A$, i.e. 
$X_{ij} = A_{ij}(e_i e_i^\tran + e_j e_j^\tran)$ where $(e_i)$ denotes the canonical basis in $\R^n$.
Then apply Corollary~\ref{cor: matrix Bernstein expectation} to the sum of mean zero matrices 
$X_{ij} - p$.    It is quick to check that $\s^2 \le pn$ and obviously $K \le 2$, and so we conclude that
\begin{equation}		\label{eq: conc from Bernstein}
\E \|A - \E A\| \lesssim \sqrt{d \log n} + \log n,
\end{equation}
where $d = np$ is the expected degree. The same argument applies more generally to 
inhomogeneous random graphs $G(n,(p_{ij}))$, and it still gives \eqref{eq: conc from Bernstein} 
when 
$$
d = \max_i \sum_j p_{ij}
$$
is the maximal expected degree. 

The logarithmic factors in bound \eqref{eq: conc from Bernstein} are not optimal, 
and can be improved by applying the result of Bandeira and van Handel (Theorem \ref{Bandeira}) to the centered adjacency matrix $A - \E A$ of an inhomogeneous random graph $G(n,(p_{ij}))$. In this case, $\s_{ij}^2 = p_{ij}$ and $K_{ij} \le 1$, so we 
obtain the following sharpening of \eqref{eq: conc from Bernstein}.

\begin{proposition}[Concentration of inhomogeneous random graphs]
  Let $A$ be the adjacency matrix of an inhomogeneous random graph $G(n,(p_{ij}))$.   Then 
  \begin{equation}		\label{eq: conc inhomogeneous}
  \E \|A - \E A\| \lesssim \sqrt{d} + \sqrt{\log n},
  \end{equation}
  where $d = \max_i \sum_j p_{ij}$ is the expected maximal degree.
\end{proposition}
In particular, if the graph is not too sparse, namely $d \gtrsim \log n$, then the optimal concentration
\eqref{eq: optimal concentration} holds, i.e. 
$$
\E \|A - \E A\| \lesssim \sqrt{d}.
$$
This recovers a result of Feige-Ofek \cite{FeiOfe05}.

A similar bound can be alternatively proved using the general result of Seginer (Theorem \ref{thm:  Seginer}).    If $A$ is the adjacency matrix of $G(n,p)$, it is easy to check that 
$\E \max_i \|A_i\|_2 \lesssim \sqrt{d} + \sqrt{\log n}$. Thus, Seginer's theorem implies the 
optimal concentration bound \eqref{eq: conc inhomogeneous} as well.
Using simple convexity arguments, one can extend this to inhomogeneous random graphs $G(n,(p_{ij}))$, 
and get the bound \eqref{eq: conc inhomogeneous} for $d = \max_{ij} np_{ij}$,
see \cite[Section~4.1]{Hajek&Wu&Xu2016}. 

One may wonder if Seginer's theorem holds for matrices with independent but not identically distributed entries.  
Unfortunately, this is not the case in general; a simple counterexample was found by Seginer \cite{Seginer2000}, see \cite[Remark~4.8]{Bandeira&Handel2016}. 
Nevertheless, it is an open conjecture of Latala that Seginer's theorem does hold if $M$ has independent {\em Gaussian} entries, see the papers \cite{Riemer&Schutt2013, vanHandel2017} and the survey \cite{vanHandelstructured2017}.

\section{Sparse networks concentrate after regularization}
\label{sec: regularization}

\subsection{Sparse networks do not concentrate} 
In the sparse regime $d=np \ll \log n$, the Bai-Yin's law for $G(n,p)$ fails. 
This is because in this case,  degrees of some vertices are much higher than the expected degree $d$.    This causes some rows of the adjacency matrix $A$ to have Euclidean norms much larger than $\sqrt{d}$, 
which in turn gives 
$$
\|A - \E A\| \gg \sqrt{d}.
$$
In other words, concentration fails for very sparse graphs; there exist outlying eigenvalues 
that escape the interval $[-2,2]$ where the spectrum of denser graphs lies according to \eqref{eq: 2}. 
For precise description of this phenomenon, see the original paper \cite{Krivelevich&Sudakov2003}, a discussion in \cite[Section~4]{Bandeira&Handel2016} and the very recent work \cite{Benaych-georgeseigenvalues2017}.


\subsection{Sparse networks concentrate after regularization}
One way to regularize a random network in the sparse regime is to remove high degree vertices altogether from the network. Indeed, \cite{FeiOfe05} showed that for $G(n,p)$, if we drop all vertices with degrees, say, larger than $2d$, then the remaining part of the network satisfies
$\|A - \E A\| =O(\sqrt{d})$ with high probability. The argument in \cite{FeiOfe05} is based on the method developed by \cite{Friedman&Kahn&Szemeredi1989} and it is extended to the IERM in \cite{Lei&Rinaldo2015,Chin&Rao&Vu2015}.

Although removal of high degree vertices restores concentration, in practice this is a bad idea, since the loss of edges associated with ``hub'' nodes in an already sparse network leads to a considerable loss of information, and in particular community detection tends to break down.   A more gentle regularization proposed in \cite{Le&Levina&Vershynin2017} does not remove high degree vertices, but reduces the weights of their edges just enough to keep the degrees bounded by $O(d)$.

\begin{theorem}[Concentration of regularized adjacency matrices]  \label{thm: main informal}

  Consider a random graph from the inhomogeneous Erd\"os-R\'enyi model $G(n,(p_{ij}))$,
  and let $d = \max_{ij} n p_{ij}$.
 
  Consider any subset of at most $10n/d$ vertices, and
  reduce the weights of the edges incident to those vertices in an arbitrary way, but so that all degrees
  of the new (weighted) network become bounded by $2d$.
 For any $r \ge 1$, with probability at least $1-n^{-r}$ the adjacency matrix $A'$ of the new weighted graph satisfies
  $$
  \|A' - \E A\| \le C r^{3/2} \sqrt{d}.
  $$
\end{theorem}

Proving concentration for this kind of general regularization requires different tools. One key result we state next is the Grothendieck-Pietsch factorization, a 
 general and well-known result in functional analysis \cite{Pietsch1978, Pisier1986,Tomczak-Jaegermann1989,Pisier2012} which 
has already been used in a similar probabilistic context \cite[Proposition 15.11]{Ledoux&Talagrand1991}.   It compares two matrix norms,
the spectral norm $\ell_2 \to \ell_2$ 
and the $\ell_\infty \to \ell_2$ norm. 

\begin{theorem}[Grothendieck-Pietsch factorization]  \label{thm:G-P}
Let $B$ be a $k \times m$ real matrix.   Then there exist positive weights $\mu_j$ with $\sum_{j = 1}^m \mu_j =1$ such that
  \begin{equation*}         \label{eq: GP factorization}
  \| B \|_{\infty\rightarrow 2} \le \| B D_\mu^{-1/2}\| \le 2 \| B \|_{\infty\rightarrow 2},
  \end{equation*}
  where $D_\mu = \diag(\mu_j)$ denotes the $m\times m$ diagonal matrix with weights $\mu_j$ on the diagonal.
\end{theorem}

\subsubsection*{Idea of the proof of Theorem~\ref{thm: main informal} by network decomposition}
\begin{figure}[htp]			
  \centering
  \begin{subfigure}[b]{0.3\textwidth}
    \includegraphics[trim=300 58 290 40,clip,width=1\textwidth]{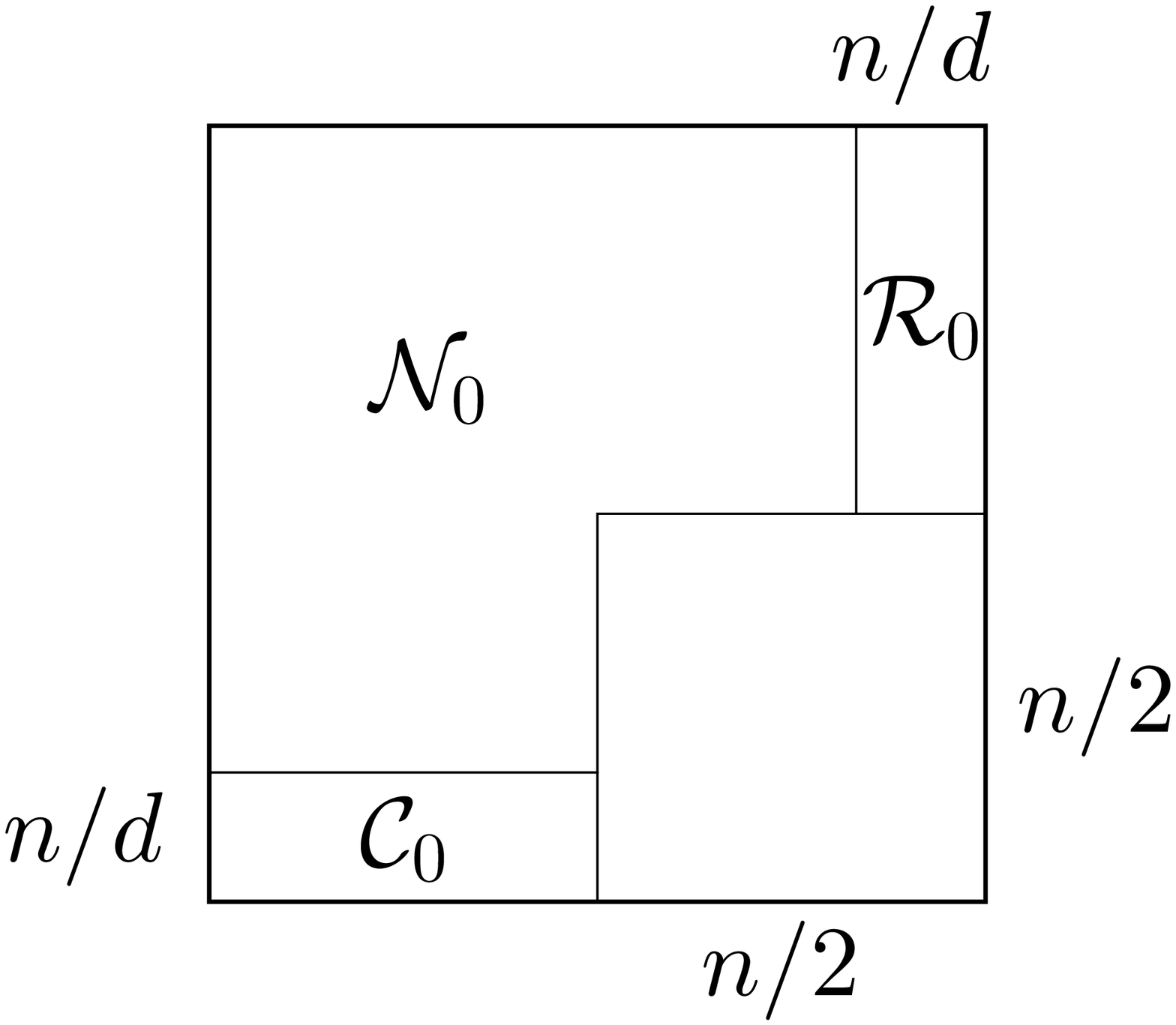}
    \caption{First step}
    \label{fig: first decomposition}
  \end{subfigure}
  \quad
  \begin{subfigure}[b]{0.3\textwidth} \qquad
    \raisebox{10pt}{\includegraphics[trim=300 100 290 40,clip,width=1\textwidth]{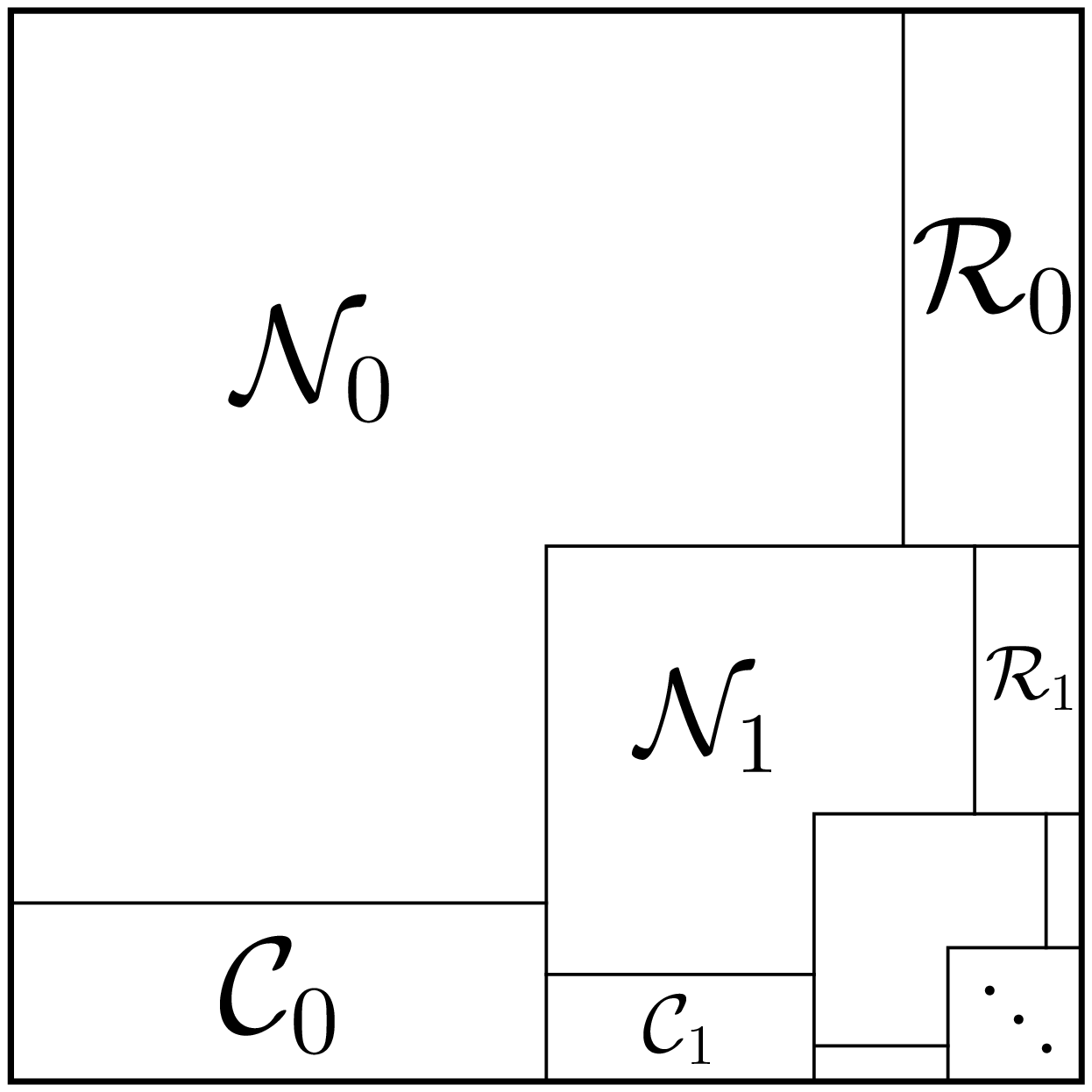}}
    \caption{Iterations}
    \label{fig: decomposition-iteration}
  \end{subfigure}
    \quad
  \begin{subfigure}[b]{0.3\textwidth} \qquad
    \raisebox{14pt}{\includegraphics[trim=300 50 290 40,clip,width=.736\textwidth]{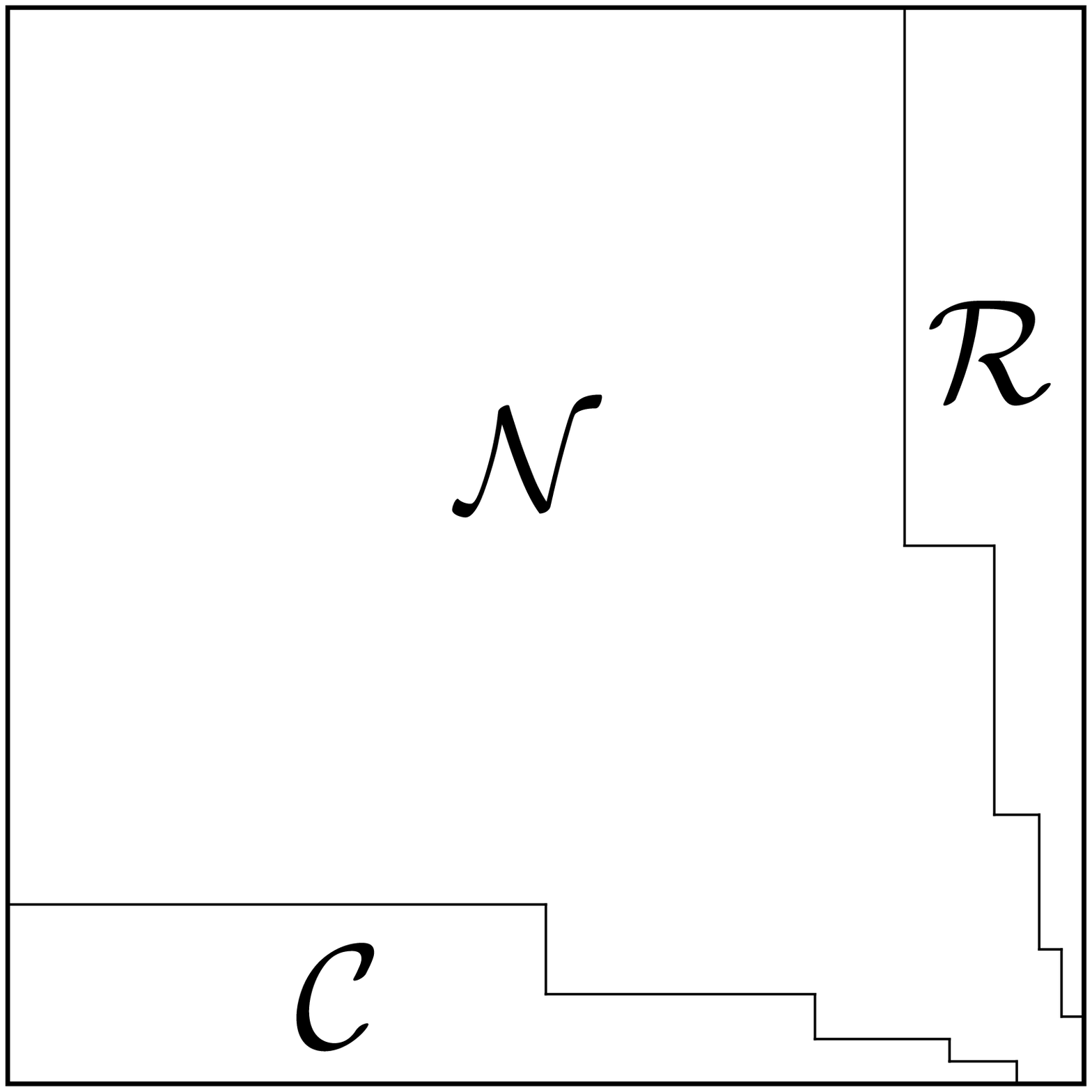}}
    \caption{Final decomposition}
    \label{fig: decomposition}
  \end{subfigure}
  \caption{Constructing network decomposition iteratively.}
\end{figure}

The argument in \cite{FeiOfe05} becomes very complicated for handling the general regularization in Theorem~\ref{thm: main informal}.    A simpler alternative approach was developed by \cite{Le&Levina&Vershynin2017} for proving Theorem~\ref{thm: main informal}. The main idea is to decompose the set of entries $[n]\times[n]$ into different subsets with desirable properties.   There exists a partition (see Figure~\ref{fig: decomposition} for illustration)   
$$
[n]\times[n] = \mathcal{N} \cup \mathcal{R} \cup \mathcal{C}
$$
such that $A$ concentrates on $\mathcal{N}$ even without regularization, while restrictions of $A$ onto $\mathcal{R}$ and $\mathcal{C}$ have small row and column sums, respectively. It is easy to see that the degree regularization does not destroy the properties of $\mathcal{N}$, $\mathcal{R}$ and $\mathcal{C}$. Moreover, it creates a new property, allowing for controlling the columns of $\mathcal{R}$ and rows of $\mathcal{C}$. Together with the triangle inequality, this implies the concentration of the entire network.

The network decomposition is constructed by an iterative procedure. We first establish concentration of $A$ in $\ell_\infty\rightarrow\ell_2$ norm using standard probability techniques. Next, we upgrade this to concentration in the spectral norm $\|(A-\E A)_{\mathcal{N}_0}\| = O(\sqrt{d})$ on an appropriate (large) subset $\mathcal{N}_0\subseteq [n]\times[n]$ using the Grothendieck-Pietsch
factorization (Theorem \ref{thm:G-P}). It remains to control $A$ on the complement of $\mathcal{N}_0$. That set is small; it can be described as a union of a block $\mathcal{C}_0$ with a small number of rows, a block $\mathcal{R}_0$ with a small number of columns and an exceptional (small) block (see Figure~\ref{fig: first decomposition}).
Now we repeat the process for the exceptional block, decomposing it into $\mathcal{N}_1$, $\mathcal{R}_1$, and $\mathcal{C}_1$, and so on, as shown in Figure~\ref{fig: decomposition-iteration}. At the end, we set $\mathcal{N} = \cup_i \mathcal{N}_i$, $\mathcal{R} = \cup_i \mathcal{R}_i$ and $\mathcal{C} = \cup_i \mathcal{C}_i$. 
The cumulative error from this iterative procedure can be controlled appropriately; see  \cite{Le&Levina&Vershynin2017} for details.  


\subsection{Concentration of the graph Laplacian}
So far, we have looked at random graphs through the lens of their adjacency matrices.   Another matrix that captures the structure of a random graph is the Laplacian.   There are several ways to define the Laplacian;  we focus on the symmetric, normalized Laplacian, $$\LL(A)=D^{-1/2} A D^{-1/2}.$$
Here $D = \diag(d_i)$ is the diagonal matrix with degrees $d_i = \sum_{j=1}^n A_{ij}$ on the diagonal.
The reader is referred to \cite{ChungFan1997} for an introduction to graph Laplacians
and their role in spectral graph theory. Here we mention just two basic facts:
the spectrum of $\LL(A)$
is a subset of $[-1,1]$, and the largest eigenvalue is always one.

In the networks literature in particular, community detection has been mainly done through spectral clustering on the Laplacian, not on the adjacency matrix.  We will discuss this in more detail in Section \ref{sec:community}, but the primary reason for this is degree normalization:  as discussed in Section \ref{sec:models}, real networks rarely have the Poisson or mixture of Poissons degree distribution that characterizes the stochastic block model;  instead, ``hubs'', or high degree vertices, are common, and they tend to break down spectral clustering on the adjacency matrix itself.   

Concentration of Laplacians of random graphs has  been studied by 
\cite{Oliveira2010,Chaudhuri&Chung&Tsiatas2012,Qin&Rohe2013,Joseph&Yu2013,Gao&Ma&Zhang&Zhou2015}.
Just like the adjacency matrix, the Laplacian is known to concentrate in the dense regime $d = \Omega(\log n)$,
and it fails to concentrate in the sparse regime. However, the reasons it fails to concentrate are different.  
For the adjacency matrix,  as we discussed, concentration fails  in the sparse case because of  high degree vertices. For the Laplacian, it is the low degree vertices that destroy concentration.   In fact, it is easy to check that when $d = o(\log n)$, the probability of isolated vertices is non-vanishing;  and each isolated vertex contributes an eigenvalue of 0 to the spectrum of $\LL(A)$, which is easily seen to
destroy concentration.

Multiple ways to regularize the Laplacian in order to deal with the low degree vertices have been proposed.  Perhaps the two most common ones are adding a small constant to all the degrees on the diagonal of $D$ \cite{Chaudhuri&Chung&Tsiatas2012}, and adding a small constant to all the entries of $A$ before computing the Laplacian.    Here we focus on the latter regularization, proposed by 
\cite{amini2013pseudo} and analyzed by \cite{Joseph&Yu2013,Gao&Ma&Zhang&Zhou2015}.
Choose $\tau > 0$ and add the same number $\tau /n$ to all
entries of the adjacency matrix $A$, thereby replacing it with
\begin{equation}
A_\tau := A + \frac{\tau}{n} \one \one^\tran
\label{eq:regLap}
\end{equation}
 Then compute the Laplacian as usual using this new adjacency matrix. 
This regularization raises all degrees $d_i$ to $d_i + \tau$, and eliminates isolated vertices, making the entire graph connected.   The original paper \cite{amini2013pseudo} suggested the choice $\tau  = \rho \bar d$, where $\bar d$ is the average node degree and $\rho \in (0,1)$ is a constant.   They showed the estimator is not particularly sensitive to $\rho$ over a fairly wide range of values away from 0 (too little regularization) and 1 (too much noise).    The choice of $\rho = 0.25$ was recommended by \cite{amini2013pseudo} but this parameter can also be successfully chosen by cross-validation on the network \cite{Li&Levina&Zhu2016}.   

The following consequence of Theorem~\ref{thm: main informal} shows
that regularization \eqref{eq:regLap}  indeed forces the Laplacian to concentrate. 
\begin{theorem}[Concentration of the regularized Laplacian]  \label{thm: Laplacian informal}
  Consider a random graph drawn from the inhomogeneous Erd\"os-R\'enyi model $G(n,(p_{ij}))$,
  and let $d = \max_{ij} n p_{ij}$.
  Choose a number $\tau>0$.
  Then, for any $r \ge 1$, with probability at least $1-e^{-r}$ we have 
  $$
  \|\LL(A_\tau) - \LL(\E A_\tau)\| \le \frac{Cr^2}{\sqrt{\tau}} \Big( 1 + \frac{d}{\tau} \Big)^{5/2}.
  $$
\end{theorem}

In the next section, we discuss why concentration of the adjacency matrix and/or its Laplacian is important in the context of community detection, the primary application of concentration in network analysis.

\section{Application to community detection}
\label{sec:community}
Concentration of random graphs has been of such interest in networks analysis primarily because it relates to the problem of community detection; see  \cite{Fortunato2010,Goldenberg2010,AbbeReview2017} for reviews of community detection algorithms and results.   We should specify that, perhaps in a slight misnomer, ``community detection'' refers to the task of assigning each node to a community (typically one and only one), not to the question of whether there are communities present, which might be a more natural use of the term ``detection''.

Most of the theoretical work linking concentration of random graphs to community detection has focused on the stochastic block model (SBM), defined in Section \ref{sec:models}, which is one of the many special cases of the general IERM we consider.
For the purpose of this paper, we focus on the simplest version of the
SBM for which the largest number of results has been obtained so far,  also known as
the balanced planted partition model $G(n,\frac{a}{n},\frac{b}{n})$.   In this model, there are $K=2$ equal-sized communities with $n/2$ nodes each.    Edges between vertices within the same community are  drawn independently with probability $a/n$, and edges between vertices in different communities with probability $b/n$.
The task is to recover the community labels of vertices from a single realization of the adjacency matrix $A$ drawn from this model.  The large literature on both the recovery algorithms and the theory establishing when a recovery is possible is very nicely summarized in the recent excellent review \cite{AbbeReview2017}, where we refer the reader for details and analogues for a general $K$ (now available for most results) and the asymmetric SBM (very few are available).  In the following subsections we give a brief summary for the symmetric $K = 2$ case which does not aim to be exhaustive.    

\subsection{Community detection phase transition} 

 {\em Weak recovery}, sometimes also called detection, means performing better than randomly guessing the labels of vertices.   The phase transition threshold for weak recovery was first conjectured in the physics literature by 
\cite{Decelle.et.al.2011}, and proved rigorously by \cite{Mossel&Neeman&Sly2014, Mossel&Neeman&Sly2014a,
  Mossel&Neeman&SlyOnConsistencyThresholds2014}, with follow-up and related work by \cite{Abbe&Bandeira&Hall2014, Massoulie:2014,Bordenave.et.al2015non-backtracking}.      The phase transition result says that there exists a polynomial time algorithm which can classify more than 50\% of the vertices correctly as $n \to \infty$ with high probability if and only if 
$$
(a-b)^2 > 2 (a+b).  
$$
Performing better than random guessing is the weakest possible guarantee of performance, which is of interest in the very sparse regime of $d = (a+b)/2 = O(1)$;  when the degree grows, weak recovery becomes trivial.   This regime has been mostly studied by physicists and probabilists;    in the statistics literature, {\em consistency} has been of more interest.   


\subsection{Consistency of community detection} 
Two types of consistency have been discussed in the literature.  Strong consistency, also known as exact recovery, means labeling {\em all} vertices correctly with high probability, which is, as the name suggests, a very strong requirement.    Weak consistency, or ``almost exact'' recovery, is the weaker and arguably more practically reasonable requirement that the fraction of misclassified vertices goes to 0 as $n \to \infty$ with high probability.   

{\em Strong consistency} was studied first, in a seminal paper \cite{Bickel&Chen2009}, as well as by \cite{Mossel&Neeman&SlyOnConsistencyThresholds2014,McS01,
Hajek&Wu&Xu2016,Cai&Li2015}.    Strong consistency is achievable, and achievable in polynomial time, if 
$$ \left| \sqrt{\frac{a}{\log n}} - \sqrt{\frac{b}{\log n}} \right| > \sqrt{2}$$
and not possible if $\left| \sqrt{a/n} - \sqrt{b/n} \right| < \sqrt{2}$.
In particular, strong consistency is normally only considered in the semi-dense regime of $d / \log n \rightarrow \infty$.   

{\em Weak consistency}, as one would expect, requires a stronger condition than weak recovery  but a weaker one than strong consistency.   Weak consistency is achievable if and only if
$$
\frac{(a-b)^2}{a+b}  = \omega (1)
$$
see for example \cite{Mossel&Neeman&SlyOnConsistencyThresholds2014}.   In particular, weak consistency is achievable in the semi-sparse regime of $d \rightarrow \infty$.   

{\em Partial recovery}, finally, refers to the situation where the fraction of misclassified vertices does not go to 0, but remains bounded by a constant below 0.5.  More specifically, partial recovery means that for a fixed $\e>0$ one can recover communities up to $\e n$ mislabeled vertices.    For the balanced symmetric case, this is true as long as
$$ \frac{(a-b)^2}{a+b} = O(1)$$
which is primarily relevant when $d = O(1)$.    Several types of algorithms are known to succeed at partial recovery in this very sparse regime, including non-backtracking walks
\cite{Mossel&Neeman&Sly2014,Massoulie:2014,Bordenave.et.al2015non-backtracking},
spectral methods \cite{Chin&Rao&Vu2015} and methods based on semidefinite programming
\cite{Guedon&Vershynin2014,Montanari&Sen2016}.   

\subsection{Concentration implies recovery}
As an example application of the new concentration results, we demonstrate how to show that  {\em regularized spectral clustering} \cite{amini2013pseudo,Joseph&Yu2013},
one of the simplest and most popular algorithms for community detection, can
recover communities in the sparse regime of constant degrees.  In general, spectral
clustering works by computing the leading eigenvectors of either the
adjacency matrix or the Laplacian, or their regularized versions, and
running the $k$-means clustering algorithm on the rows of the $n \times k$  matrix of leading eigenvectors to
recover the node labels.    In the simplest case of the balanced $K =2$ model
$G(n,\frac{a}{n},\frac{b}{n})$, one can simply assign nodes to two 
communities according to the sign of the entries of the
eigenvector $v_2(A')$ corresponding to the second smallest
eigenvalue of the (regularized) adjacency matrix $A'$.

Let us briefly explain how concentration results validate recovery from the regularized adjacency matrix or regularized Laplacian.   
If concentration holds and the regularized matrix $A'$ is shown to be close to $\E A$, then standard perturbation theory (i.e., the Davis-Kahan theorem, see e.g. \cite{Bhatia1996}) implies that $v_2(A')$ is close to $v_2(\E A)$, and in particular, the signs of these two eigenvectors must agree on most vertices.
An easy calculation shows that the signs of $v_2(\E A)$ recover the communities exactly: 
the eigenvector corresponding to the second smallest eigenvalue of $\E A$ (or the second largest of $\LL(A)$) is a positive constant on one community and a negative constant on the
other. Therefore, the signs of $v_2(A')$ recover communities up to a small fraction of misclassified vertices and, as always, up to a permutation of community labels. This argument remains valid if we replace the regularized adjacency matrix $A'$ with regularized Laplacian $\LL(A_\tau)$.

\begin{corollary}[Partial recovery from a regularized adjacency matrix for sparse graphs]\label{cor: detectionA}
Let $\e>0$ and $r\ge 1$. Let $A$ be the adjacency matrix drawn from the stochastic block model $G(n,\frac{a}{n},\frac{b}{n})$. Assume that
\begin{equation*}
  (a-b)^2>C(a+b)
\end{equation*}
where $C$ is a constant depending only on $\varepsilon$ and $r$.
For all nodes with degrees larger than $2a$,
reduce the weights of the edges incident to them in an arbitrary way, but so that all degrees
of the new (weighted) network become bounded by $2a$, resulting in a new matrix $A'$.   
Then with probability at least $1-e^{-r}$, the signs of the entries of the eigenvector corresponding to the second smallest eigenvalue of $A'$ correctly estimate the partition into two communities, up to at most $\e n$ misclassified vertices.  
\end{corollary}

\begin{corollary}[Partial recovery from a regularized Laplacian for sparse graphs] \label{cor: detectionL}
Let $\e>0$ and $r\ge 1$. Let $A$ be the adjacency matrix drawn from the stochastic block model $G(n,\frac{a}{n},\frac{b}{n})$. Assume that
\begin{equation}\label{eq: ab}
  (a-b)^2>C(a+b)
\end{equation}
where $C$ is a constant depending only on $\varepsilon$ and $r$.
Choose $\tau$ to be the average degree of the graph, i.e. $\tau=(d_1+\cdots+d_n)/n$.
Then with probability at least $1-e^{-r}$, 
 the signs of the entries of the eigenvector corresponding to the second largest eigenvalue of $\LL(A_\tau)$ correctly estimate the partition into the two communities, up to at most $\e n$ misclassified vertices.
\end{corollary}

\begin{figure}[!ht]
  \centering
  \includegraphics[trim=5 5 20 10,clip,width=1\textwidth]{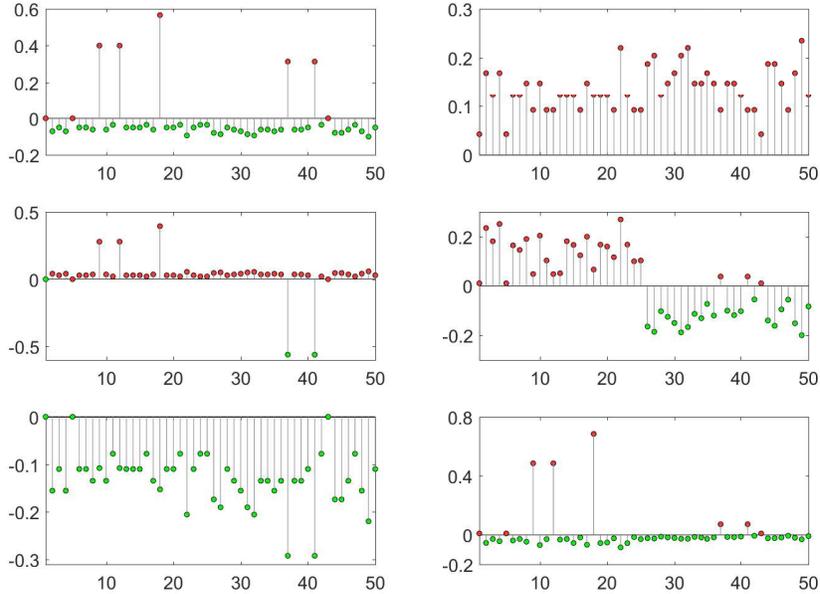}\\
  \caption{Three leading eigenvectors (from top to bottom) of the Laplacian (left) and the regularized Laplacian (right).   The network is generated from $G(n,\frac{a}{n},\frac{b}{n})$ with $n=50$, $a=5$ and $b=0.1$.  
Nodes are labeled so that the first 25 nodes belong to one community and the rest to the other community.    Regularized Laplacian is computed from $A + 0.1 \bar d / n \one \one^\tran$.   
}
\label{fig: eigenvectors}
\end{figure}

As we have discussed,  the Laplacian is typically preferred over the adjacency matrix in practice, because the variation in node degrees is reduced by the normalization factor $D^{-1/2}$ \cite{Sarkar&Bickel2015}.   Figure~\ref{fig: eigenvectors} shows the effect of regularization for the Laplacian of a random network generated from $G(n,\frac{a}{n},\frac{b}{n})$ with $n = 50$, $a=5$ and $b=0.1$. For plotting purposes, we order the nodes so that the first $n/2$ nodes belong to one community and the rest belong to the other community.  Without regularization, the two leading eigenvectors of the Laplacian localize around a few low degree nodes, and therefore do not contain any information about the global community structure. In contrast, the second leading eigenvector of the regularized Laplacian (with $\tau = 0.1 \bar d$) clearly reflects the communities, and the signs of this eigenvector alone recover community labels correctly for all but three nodes.

\section{Discussion}

 Great progress has been made in recent years, and yet many problems remain open.   Open questions on community detection under the SBM, in terms of exact and partial recovery and efficient (polynomial time) algorithms are discussed in \cite{AbbeReview2017}, and likely by the time this paper comes out in print, some of them will have been solved.    Yet the focus on the SBM is unsatisfactory for many practitioners, since not many real networks fit this model well.   Some of the more general models we discussed in Section~\ref{sec:models} fix some of the problems of the SBM, allowing for heterogeneous degree distributions and overlapping communities, for instance.   A bigger problem lies in the fixed $K$ regime; it is not realistic to  assume that as the size of the network grows, the number of communities remains fixed.   A more realistic model is the ``small world'' scenario, where the size of communities remains bounded or grows very slowly with the number of nodes, the number of communities grows, and connections between many smaller communities happen primarily through hub nodes.   Some consistency results have been obtained for a growing $K$, but we are not aware of any results in the sparse constant degree regime so far.    An even bigger problem is presented by the so far nearly universal assumption of independent edges;  this assumption violates commonly observed transitivity of friendships (if A is friends with B and B is friends with C, A is more likely to be friends with C).   There are other types of network models that do not rely on this assumption, but hardly any random matrix results apply there.   Ultimately, network analysis involves a lot more than community detection:   link prediction, network denoising, predicting outcomes on networks, dynamic network modeling over time, and so on.   We are a long way away from establishing rigorous theoretical guarantees for any of these problems to the extent that we have for community detection, but given how rapid progress in the latter area has been, we are hopeful that  continued interest from the random matrix community will help shed light on other problems in network analysis.   

\bibliography{allref}
\bibliographystyle{abbrv}

\end{document}